\documentclass[12pt]{amsart}

\usepackage{amsmath, amssymb, graphics}

\newtheorem{theorem}{Theorem}[section]

\begin{document}

\title{On the hyperbolic metric of certain domains}

\author[Aimo Hinkkanen]{Aimo Hinkkanen}
\address{Aimo Hinkkanen: Department of Mathematics, 
University of Illinois at Urbana--Champaign, U.S.A.}
 \email{aimo@illinois.edu}
\author[Matti Vuorinen]{Matti Vuorinen}
\address{Matti Vuorinen: Department of Mathematics and Statistics, University of Turku, FI-20014 Turku, Finland}
\email{vuorinen@utu.fi}

\dedicatory{Dedicated to the memory of Peter L.\ Duren}

\thanks{}


\abstract We prove that if $E$ is a compact subset of the unit disk ${\mathbb D}$ in the complex plane, if $E$ contains a sequence of distinct  points $a_n\not= 0$ for $n\geq 1$ such that $\lim_{n\to\infty} a_n=0$ and  for all $n$ we have
$
|a_{n+1}| \geq \frac{1}{2} |a_n|  
$, and if
$G={\mathbb D} \setminus E$ is connected and   $0\in \partial G$, then  there is a constant $c>0$ such that for all $z\in G$ we have
$
\lambda_{G } (z) \geq c/|z|
$
where $\lambda_{G } (z)$ is the density of the hyperbolic metric in $G$. 
\endabstract

\maketitle

\section{Introduction}

We denote by ${\mathbb D} = \{ z\in {\mathbb C} \colon |z|<1\}$ 
the unit disk in the complex plane ${\mathbb C}$. 
For $z\in {\mathbb C}$ and $r>0$, we write 
$B(z,r) = \{ w\in {\mathbb C} \colon |w-z|<r \}$ and 
$S(z,r) = \{ w\in {\mathbb C} \colon |w-z|=r \}$. 
If $G$ is a domain in ${\mathbb C}$ with a hyperbolic metric, 
we denote by $\lambda_{G } (z)$, for $z\in G$, the 
density of the hyperbolic metric in $G$. We use the normalization under which
\begin{equation} \label{c0}
\lambda_{{\mathbb D} } (z) = \frac{1}{1-|z|^2}
\end{equation}
for all $z\in {\mathbb D}$. 

The hyperbolic metric has numerous applications to geometric function theory. Concrete estimates for the metric or for its density function are useful, e.g., in the study of distortion under analytic functions and quasiregular mappings \cite{BP}. We prove the following theorem which improves the well-known a priori estimate of the type $d/(|z| \log |z|)$, for some constant $d>0$, valid by the comparison principle (see, e.g., \cite{A}, p.~17).

\begin{theorem} \label{th1}
Let $E$ be a compact subset of ${\mathbb D}$ and suppose that $E$ contains a sequence of distinct   points $a_n\not= 0$ for $n\geq 1$ such that $\lim_{n\to\infty} a_n=0$ and such that  for all $n$ we have
\begin{equation} \label{c1}
|a_{n+1}| \geq \frac{1}{2} |a_n|  .
\end{equation}
Suppose that $G={\mathbb D} \setminus E$ is connected, hence a domain, and that $0\in \partial G$. Then  there is a constant $c>0$ such that for all $z\in G$ we have
\begin{equation} \label{c2}
\lambda_{G } (z) \geq \frac{c}{|z|}   .
\end{equation}
We may take
\begin{equation} \label{c2a}
c = \min \left\{  \frac{1}{  2\sqrt{2} \left(\kappa + \log \frac{4}{\delta} \right) },  \frac{1}{  2\sqrt{2} \left(\kappa + 5 \log 2  \right) }  \right\}   ,
\end{equation}
where
\begin{equation} \label{c3}
\delta = \max \{ |a_n| \colon n\geq 1\}  
\end{equation}
and
\begin{equation} \label{c72}
\kappa = 4 + \log (3 + 2 \sqrt{2}) \approx 5.7627...  .
\end{equation}
\end{theorem}

We remark that
\begin{equation} \label{c72a}
\frac{1}{  2\sqrt{2} \left(\kappa + 5 \log 2  \right) } \approx 0.03831...  
\end{equation}
and that
$$
\frac{1}{  2\sqrt{2} \left(\kappa + \log \frac{4}{\delta} \right) } <  \frac{1}{  2\sqrt{2} \left(\kappa + 5 \log 2  \right) }
$$
when $0<\delta<1/8$.

In Section~\ref{ex} we give an example to show that the constant $c$ necessarily depends on $E$. More precisely, in simple terms, the correct order of magnitude for $c$ is
$c \sim C_1/( \log \frac{1}{ \delta' }   ) $ where $C_1\in (0,1)$ is an absolute constant and $\delta' = \min\{1/2,  \max \{ |a_n| \colon n\geq 1\} \} \in (0,1) $. 

The assumptions on the sequence $a_n$ imply that $0\in E$ since $E$ is closed. Without further assumptions on $E$, other than that $E$ contains the sequence $a_n$,  there is no need for ${\mathbb D} \setminus E$ to be connected, and if it is not connected, then  ${\mathbb D} \setminus E$ could have components $\Omega$ such that $\partial \Omega \cap \partial {\mathbb D} = \emptyset$, and in addition $\Omega$ might be multiply connected, even of infinite connectivity. For this reason we assume that $G={\mathbb D} \setminus E$ is connected, hence a domain, possibly multiply connected.

If there is $r>0$ such that $E$ contains the disk $B(0,r)$, then $c/|z|\leq c/r$ for all $z\in G$ while $\lambda_{G } (z)$ has a positive minimum on $G$. Hence it is trivial that (\ref{c2}) holds, even though $c$ may depend on $G$ in a manner that is not specified here. To avoid this situation, we have assumed that $0\in \partial G$. Then for every $r\in (0,1)$, the domain $G$ intersects the circle $S(0,r)$.

As an application of Theorem~\ref{th1} we obtain the following result, which is similar to Littlewood's distortion theorem (\cite{L}, Theorem~24, p.~502), see also \cite{JV}, Theorem~3.1, p.~520.

\begin{theorem} \label{th2}
Let $G$ be a domain satisfying the assumptions of Theorem~\ref{th1} and let $f$ be a non-constant analytic function of  ${\mathbb D}$ into $G$. Then for all $z\in {\mathbb D}$ we have
\begin{equation} \label{d1}
c \left| \log \frac{    |f(0 )|   } {     |f(z )|    }     \right| \leq 
\rho_G(f(0),f(z)) \leq \rho_{   {\mathbb D}   } (0,z) = \frac{1}{2} \log \frac  { 1+|z| }  { 1 - |z| } 
\end{equation}
where $\rho_G$ stands for the hyperbolic distance defined by $\lambda_G$ and $c$ is the constant of Theorem~\ref{th1}. In particular, for all $z\in {\mathbb D}$ we have
\begin{equation} \label{d2}
|f(z)| \geq |f(0)| \left(   \frac  { 1 - |z| }  { 1 + |z| }       \right)^{1/(2c)}  .   
\end{equation}
\end{theorem}

{\bf Proof of Theorem~\ref{th2}.} The first inequality (\ref{d1}) follows by the definition of the hyperbolic distance together with (\ref{c2}),  and the second inequality in (\ref{d1}) is a basic property of analytic functions. \qed

We will make extensive use of the following estimate for the hyperbolic metric, due to Beardon and Pommerenke (\cite{BP}, Theorem~1, p.~477). 

{\bf Theorem A.} {\sl 
Let $D$ be a plane domain with a hyperbolic metric. For $z\in D$, define
\begin{equation} \label{c71}
L=L(z)= \inf \left\{ \left|  \log \left|  \frac{ z-a} { a-b  }   \right| \right| \colon 
a,b\in \partial D, \,\,  |z-a| = d(z,\partial D)  \right\}
\end{equation}
where $d(z,S)$ is the Euclidean distance from $z$ to the set $S$. Then
\begin{equation} \label{c7}
\frac{  1} {  2 \sqrt{2} d(z,\partial D) (\kappa +L)  } \leq \lambda_{D }(z)
 \leq \frac{\kappa +  (\pi/4)    } { d(z,\partial D) (\kappa + L)      }   
\end{equation} 
where $\kappa$ is given by (\ref{c72}).  
}

\section{Proof of Theorem~\ref{th1}}

\subsection{Preliminary observations.}
We first make some observations. Since $\lim_{n\to\infty} a_n=0$, the quantity
$\delta$ given by (\ref{c3})
is well defined, and there is $j\geq 1$ such that $|a_j|=\delta$, so $0<\delta<1$.   

We next claim that for each $n\geq 0$, there is $k$ such that
\begin{equation} \label{c4}
2^{-(n+1)} \delta < |a_k| \leq   2^{-n} \delta  .  
\end{equation}
We prove this by induction on $n$. When $n=0$, the claim is valid since $|a_j|=\delta$. Suppose that the claim is valid for a certain $n\geq 0$. Thus there is at least one $m$ such that $2^{-(n+1)} \delta < |a_m| \leq   2^{-n} \delta$. 
There are only finitely many $k$ such that $2^{-(n+1)} \delta < |a_k| $. Among these $k$, let $k_0$ be maximal. By the definition of $k_0$, we then have $|a_k|\leq 2^{-(n+1)} \delta$ for all $k>k_0$, and in particular,
\begin{equation} \label{c5}
|a_{k_0+1}| \leq 2^{-(n+1)} \delta  .  
\end{equation}
On the other hand, since $|a_{p+1}| \geq \frac{1}{2} |a_p| $ for all $p$ and since $2^{-(n+1)} \delta < |a_{k_0}| $, we have
$|a_{k_0+1}| \geq (1/2)|a_{k_0}| > 2^{-(n+2)} \delta$. Putting these results together we see that
$$
2^{-(n+2)} \delta < |a_{k_0+1}| \leq 2^{-(n+1)} \delta  .
$$
This completes the induction proof.

\subsection{Closest boundary point on the unit circle.}
We apply the left hand inequality (\ref{c7}) taking $D$ to be $G$.

Fix $z\in G$, so $z\not= 0$. Consider the points $\zeta\in \partial G$ such that
\begin{equation} \label{cb}
|z-\zeta| = d(z,\partial G)  .
\end{equation} 
For the rest of the proof, $\zeta$ will always denote   a point satisfying (\ref{cb}). 
Since $0\in \partial G$, we always have $d(z,\partial G) \leq |z|$, hence $|z-\zeta|  \leq |z|$.
There may be more than one such point $\zeta$. Suppose first that we may choose $\zeta$ so that $|\zeta|=1$. Then $|z-\zeta|  = 1-|z|$, for otherwise the line segment joining $z$ to $z/|z|$ intersects $\partial G$. Then also $|z|\geq 1/2$. For if $|z|<1/2$, then, since $0\in E$, there would be $\zeta$ with $|z-\zeta|\leq |z|<1/2<1-|z|$, a contradiction. There is $w\in \partial {\mathbb D} \subset \partial G$ with $|w-\zeta|=|z-\zeta|$. Choosing $a=\zeta$ and $b=w$ in (\ref{c71}), we see that $L(z)=0$. 
In this case $|z-\zeta| =1-|z|\leq |z|$, so that by (\ref{c7}), 
$$
\lambda_{G }(z) \geq \frac{  1} {  2 \sqrt{2} \kappa (1-|z|)   }\geq \frac{  1} {  2 \sqrt{2} \kappa |z|   } ,
$$
which gives (\ref{c2}) with 
\begin{equation} \label{c8a}
c=1/(2 \sqrt{2} \kappa)  .  
\end{equation}  

\subsection{Closest boundary point inside the unit disk.}
Suppose then that we cannot choose $\zeta$ in (\ref{cb}) so that $|\zeta|=1$. 

Consider what we need to prove. By (\ref{c7}), the inequality (\ref{c2}) will hold provided that
\begin{equation} \label{c8}
\frac{  1} {  2 \sqrt{2} |z-\zeta| (\kappa +L)  }\geq \frac{c}{|z|},
\end{equation} 
that is,
\begin{equation} \label{c9}
 2 \sqrt{2} c (\kappa +L)  \leq \frac {|z|} { |z-\zeta|  } .
\end{equation} 
Since $|z-\zeta| \leq |z|$, the right hand side is $\geq 1$, which will be a useful fact. Considering the definition of $L$, we see that (\ref{c9}) holds if  we can choose, apart from $\zeta$, a point $b\in \partial G$ such that 
\begin{equation} \label{c9a}
 2 \sqrt{2} c \left(\kappa +\left|  \log \left|  \frac{ z-\zeta} { \zeta-b  }   \right| \right|  \right)  \leq \frac {|z|} { |z-\zeta|  } .
\end{equation} 

No matter how we choose $\zeta\in \partial G$ satisfying (\ref{cb}), we can always find $b\in \partial G$ such that $|\zeta-b|\geq \delta/2$. For if $|\zeta|\geq \delta/2$, we can choose $b=0$. If $|\zeta| < \delta/2$, we can choose $b$ as follows. Recall that $ |a_j | = \delta$. On the line segment from $a_j$ to $a_j/|a_j|$, starting from $a_j$, we find the first point, possibly $a_j$ itself, that lies on $\partial G$. We take this point to be $b$. Then $|b|\geq \delta$ so that $|\zeta-b|\geq |b|-|\zeta|>\delta/2$. We will make this choice for $b$ in some, but not all, of the situations considered below.

\subsection{The case $|z-\zeta|\geq \delta/2$.}
Suppose first that (\ref{cb}) holds and $|z-\zeta|\geq \delta/2$. Then we choose $b$ as above, so that each of the quantities $|z-\zeta|$ and $|b-\zeta|$ lies in the interval $[\delta/2,2)$. Then 
$$
\left|  \log \left|  \frac{ z-\zeta} { \zeta-b  }   \right| \right|  \leq  \log \frac{4} { \delta}   .
$$ 
Since $|z|/|z-\zeta| \geq 1$, it follows that we can choose 
\begin{equation} \label{c99}
c= \frac{ 1 } {  2 \sqrt{2}  \left(\kappa + \log \frac{4} { \delta}   \right)  } 
\end{equation} 
to be a positive constant depending on $\delta$ only so that (\ref{c9a}) holds. 

\subsection{The case $|z-\zeta|< \delta/2$ and $|z|\geq \delta/2$.}
Suppose then that $|z-\zeta| <  \delta/2$. Suppose first that also $|z|\geq \delta/2$. Choose $b$ as above so that $|\zeta-b|\geq \delta/2$. Then
$$
\left|  \log \left|  \frac{ z-\zeta} { \zeta-b  }   \right| \right|  = 
\log \left|  \frac{ b-\zeta} { \zeta-z  }   \right| \leq 
\log 2 - \log |z-\zeta| .
$$
Hence (\ref{c9a}) will hold provided that
\begin{equation} \label{c9b}
 2 \sqrt{2} c \left(\kappa +\log 2 - \log |z-\zeta| \right)  \leq \frac {\delta/2} { |z-\zeta|  } .
\end{equation} 
Consider $t=|z-\zeta|$ to be a variable in $(0,\delta/2)$ and define
$$
\varphi(t) = \frac {\delta/2} { t  } - 2 \sqrt{2} c \left(\kappa +\log 2 - \log t \right)  .
$$
Then
$$
\varphi'(t) = \frac {-\delta/2} { t^2  } + 2 \sqrt{2} c \frac{1}{t}  .
$$
We have $\varphi'(t) =0$ if, and only if, $t=\delta /(4\sqrt{2}c)$. If we choose $c$ to satisfy, in addition to other conditions, that $c\leq 1/( 2 \sqrt{2} )$, then $\delta /(4\sqrt{2}c)
\geq \delta/2$, so that $\varphi'(t) <0$ and hence $\varphi$ is decreasing on $(0,\delta/2)$. We have
$$
\varphi(\delta/2) = 1 - 2 \sqrt{2} c \left(\kappa +\log 4 - \log \delta \right)
$$
so that $\varphi(\delta/2) \geq 0$ provided that $c$ is small enough, depending on $\delta$ only, for example, if $c$ has the value given by (\ref{c99}).
Then we have $\varphi(t)>0$ for all $t\in (0,\delta/2)$, and hence (\ref{c9b}) (and hence also (\ref{c9a})) holds. 

\subsection{The case $|z-\zeta| <  \delta/2$ and $|z| < \delta/2$, subcase $|z-\zeta| \leq |z|/8$.}
It remains to consider the case when (\ref{cb}) holds and $|z-\zeta| <  \delta/2$ and $|z| < \delta/2$. 
Then $|\zeta|<\delta$. Now we consider two subcases depending on whether $|z-\zeta| \leq |z|/8$ or $|z-\zeta| >|z|/8$. 

Consider first the subcase when $|z-\zeta| \leq |z|/8$ so that $|z|/|z-\zeta|\geq 8$. 
There is a unique $n\geq 0$ such that $2^{-(n+1)} \delta < |\zeta| \leq   2^{-n} \delta$.  Choose $k_1$  so that  $2^{-(n+3)} \delta < |a_{k_1}| \leq   2^{-(n+2)} \delta$. Choose a point $w\in G\cap S(0,2^{-(n+2)} \delta)$. Join $w$ to $a_{k_1}$ by a curve consisting of the shorter arc of the circle $S(0,2^{-(n+2)} \delta)$ from $w$ to $|w|a_{k_1}/|a_{k_1}|$, followed by a radial line segment from $|w|a_{k_1}/|a_{k_1}|$ to $a_{k_1}$. Let $b$ be the first point of $\partial G$ that we encounter when we traverse this curve starting from $w$. Then 
$$
2^{-(n+3)} \delta < |b| \leq   2^{-(n+2)} \delta  
$$
so that
$$
|b-\zeta| \geq |\zeta| - |b| \geq 2^{-(n+1)} \delta - 2^{-(n+2)} \delta
= 2^{-(n+2)} \delta  
$$
and
$$
|b-\zeta| \leq |\zeta| + |b| \leq  2^{-n} \delta +  2^{-(n+2)} \delta
< 2^{-(n-1)} \delta < 4 |\zeta| .
$$
Furthermore,
$$
|\zeta| - |z|/8  \leq |\zeta| - |\zeta-z|   \leq |z| \leq |\zeta| + |\zeta-z| \leq |\zeta| + |z|/8
$$
so that
$$
7|z|/8  \leq |\zeta| \leq 9|z|/8 
$$
and hence also
$$
|b-\zeta| \leq 4 |\zeta| \leq 9|z|/2  .
$$
Thus we have
$$
|b-\zeta| \geq 2^{-(n+2)} \delta \geq \frac{1}{4} |\zeta| \geq \frac{7 } {32 } |z| 
\geq  \frac{ |z| } {8} \geq  |z-\zeta|   .
$$
It follows that
$$
\left|  \log \left|  \frac{ z-\zeta} { \zeta-b  }   \right| \right|  = 
\log \left|  \frac{ b-\zeta} { \zeta-z  }   \right| \leq 
\log  \frac{  9|z|/2  } { |\zeta-z|  }   \leq
\log (9/2) + \log t  ,
$$
where we have written
$$
t = \frac{ |z| } {  |\zeta-z|   } \geq 8 .
$$
Now (\ref{c9a}) will hold provided that
\begin{equation} \label{c9c}
 2 \sqrt{2} c \left(\kappa +\log (9/2) + \log t  \right)  \leq t .
\end{equation} 
Write
$$
\varphi(t) = t - 2 \sqrt{2} c \left(\kappa +\log (9/2) + \log t  \right) .
$$
Then
$$
\varphi'(t) = 1 - \frac{  2 \sqrt{2} c } {t}   .
$$
Hence $\varphi'(t) = 0$ if, and only if, $t= 2 \sqrt{2} c$. If $c<1$ (say), we have $\varphi'(t) > 0$ for all $t>8$, and hence $\varphi(t) $ is increasing on $(8,\infty)$. We have $\varphi(8) =   8 - 2 \sqrt{2} c \left(\kappa +\log (9/2) + \log 8  \right)$ so that $\varphi(8) \geq 0$ (and hence $\varphi(t) >0$ for all $t>8$) provided that 
\begin{equation} \label{c9dd}
c \leq \frac{ 2 \sqrt{2} } { \kappa + 2\log 6  } \approx  0.302626 .
\end{equation}
Thus, with such a choice for $c$, the inequality (\ref{c9a}) will hold in this case.

\subsection{The case $|z-\zeta| <  \delta/2$ and $|z| < \delta/2$, subcase $|z-\zeta| > |z|/8$.}
In the remaining subcase (\ref{cb}) holds and we have $|z-\zeta| <  \delta/2$, $|z| < \delta/2$, and $|z-\zeta| >|z|/8$ (and also $|z-\zeta| \leq |z|$ as always). 
Now 
\begin{equation} \label{c9ddd}
1\leq \frac {|z|} { |z-\zeta|  }  \leq  8  ,
\end{equation}
so that (\ref{c9a})
follows, and  is essentially equivalent to the statement that, 
 for a suitable choice of $b\in \partial G$, we have
\begin{equation} \label{c9ee}
\left|  \log \left|  \frac{ z-\zeta} { \zeta-b  }   \right| \right|  
\leq \frac{1} { 2 \sqrt{2} c }   - \kappa  .
\end{equation}
Due to (\ref{c9ddd}),  we may consider, instead of the left hand side of (\ref{c9ee}), the quantity
$$
\left|  \log \left|  \frac{ z} { \zeta-b  }   \right| \right|  
$$
and so we seek $b$ such that
\begin{equation} \label{c9d}
c_1 |z| \leq | \zeta-b | \leq c_2 |z|
\end{equation}
for some positive constants $c_1$ and $c_2$. 
Further,
$$
|\zeta| \leq |z| + |z-\zeta| \leq 2 |z| 
$$
and
$$
|z-\zeta| \geq |z|/8 \geq |\zeta|/16. 
$$
There is a unique $n\geq 1$ such that $2^{-(n+1)} \delta < |z| \leq   2^{-n} \delta$. 

If $|\zeta|\geq |z|/4$, we choose $b=0$ so that $| \zeta-b | = |\zeta| $ and note that $|z|/4\leq |\zeta|\leq 2|z|$. Hence (\ref{c9d}) holds. 

If $|\zeta| < |z|/4\leq 2^{-(n+2)} \delta$, we find $b\in \partial G$ with $2^{-(n+1)} \delta < |b| \leq   2^{-n} \delta$. (We proceed as before. We find $a_k$ such that $2^{-(n+1)} \delta < |a_k| \leq   2^{-n} \delta$. Then we join $z$ to $a_k$ by a curve  consisting of the shorter arc of the circle $S(0,|z|)$ from $z$ to $|z|a_{k}/|a_{k}|$, followed by a radial line segment from $|z|a_{k}/|a_{k}|$ to $a_{k}$. Let $b$ be the first point of $\partial G$ that we encounter when we traverse this curve starting from $z$.) Then
$$
| \zeta-b | \geq |b| - |\zeta| \geq 2^{-(n+1)} \delta - 2^{-(n+2)} \delta
= 2^{-(n+2)} \delta \geq (1/4) |z| 
$$
and
$$
| \zeta-b | \leq |b| + |\zeta| \leq 2^{-n} \delta + 2^{-(n+2)} \delta
< 2^{-(n-1)} \delta \leq 4 |z| .
$$
Thus (\ref{c9d}) holds. 

A careful consideration of the above estimates shows that we always have
$$
\left|  \log \left|  \frac{ z-\zeta} { \zeta-b  }   \right| \right|  \leq \log 32
$$
in (\ref{c9ee}), so that we may take
\begin{equation} \label{c9ef}
c = \frac{1}{  2\sqrt{2} \left(\kappa + 5 \log 2  \right) } .
\end{equation}

We have now dealt with all cases and therefore have proved (\ref{c2}). Considering the various values of $c$ in (\ref{c8a}), (\ref{c99}), (\ref{c9dd}), and (\ref{c9ef}), we find that we may choose $c$ as given by (\ref{c2a}). 
This completes the proof of Theorem~\ref{th1}. 

{\bf Remark.}  It was not necessary to consider separately the case when we may choose $\zeta$ in (\ref{cb}) so that $|\zeta|=1$. That separate consideration merely showed that in that case we may take $c$ to be an absolute constant, independent of $\delta$.

\section{An example} \label{ex}

Suppose that $0<\delta<1/4$ and let $E$ be the line segment $[0,\delta]$. Set $G={\mathbb D}\setminus E$ and consider $z=\sqrt{\delta}\in G$. We could calculate $\lambda_G(z)$ explicitly, but the following estimates will be sufficient. By the right hand inequality in (\ref{c7}), we have
\begin{equation} \label{c9e}
\lambda_G(z) \leq \frac{\kappa +  (\pi/4)    } { d(z,\partial G) (\kappa + L)      }   .
\end{equation} 

Now for our $G$ and for $z= \sqrt{\delta} < 1/2$, we have $\zeta=\delta$, $|z-\zeta|=  \sqrt{\delta} -\delta$, and (as a bit of thinking shows) $L=\log \frac{   \sqrt{\delta}-\delta } { \delta  }$. Hence
if we are to have $c/|z|\leq \lambda_G(z)$ in this case, we must have
$$
 c\leq \frac{\kappa +  (\pi/4)    } { (1- \sqrt{\delta})  \left(\kappa + \log \frac{   \sqrt{\delta}-\delta } { \delta  } \right)      } .
$$
Considering what happens when $\delta$ decreases to $0$, we see that 
we need to have
\begin{equation} \label{c5a}
c \leq \frac{ C } {  \log \frac{1}{ \delta }   }    
\end{equation}
for some absolute constant $C>0$.

Indeed, 
Theorem~\ref{th1} shows that with $\delta'=\min\{\delta,1/2\}$, say, we may always choose $c=C_1/(\log (1/\delta'))$ for some absolute constant $C_1<1$.

\end{document}